\documentclass{article}

\usepackage{amssymb}

\newtheorem{theorem}{Theorem}[section]
\newtheorem{lemma}[theorem]{Lemma}

\newcommand{\R}{\mathbf{R}}       % Nombres réels
\newcommand{\N}{\mathbf{N}}       % Nombres entiers
\newcommand{\calL}{\mathcal{L}}   % ``convergence''
\newcommand{\calN}{\mathcal{N}}   % ``bornitude''
\newcommand{\calJ}{\mathcal{J}}   % ``finitude de l'énergie''
\newcommand{\Hyp}{(\mathcal{H})}    % l'hypothèse d'uniformité
\newcommand{\un}[1]{\mathbf{1}_{(#1)}}   % Fonction caractéristique
\newcommand{\simsubset}{\stackrel{\sim}{\subset}}  % Presque inclus
\newcommand{\dS}{\partial S}      % Bord

\begin{document}

\title{Non-tangential, radial and stochastic asymptotic properties 
       of harmonic functions on trees
       \footnote{\textsl{Math. Classif.}~:
          05C05, 31C20, 31C35, 60J15, 60J50.}}
\author{Fr\'ed\'eric Mouton}

\maketitle

\begin{abstract}
  For a harmonic function on a tree with random walk whose  
  transition probabilities are 
  bounded between two constants in $(0,1/2)$, 
  it is known that the radial and stochastic
  properties of convergence, boundedness and finiteness
  of energy are all a.s. equivalent.
  We prove here that the analogous non-tangential properties 
  are a.e. equivalent to the above ones.  
\end{abstract}

\vspace{5mm}

We are interested in the comparison between some non-tangential 
asymptotic properties of harmonic functions on a tree 
and the corresponding radial properties, using analogous stochastic ones. 
We proved in a previous work~\cite{mout3}, under 
a reasonable uniformity hypothesis, the almost sure 
equivalence between different radial and stochastic properties:
convergence, boundedness and finiteness of the energy. 
The probabilistic-geometric methods, adaptated from those we used 
in the setting of
manifolds of negative curvature~\cite{mout2}, were flexible and 
presumed to extend to the non-tangential case for trees.

A recent article~\cite{atan-pica} shows by combinatorial methods the equivalence
of the three non-tangential corresponding properties in the 
particular case of homogeneous trees. 
It seems to be time to show explicitely that our methods give
in a swift way the non-tangential results for general trees satisfying the
uniformity hypothesis above.

We use our previous results to compare the non-tangential notions with 
the radial and stochastic ones: we prove on one hand that the stochastic
convergence implies the non-tangential convergence in the section~\ref{csicnt}
and on the other hand that the non-tangential boundedness implies almost
surely the finiteness of the non-tangential energy in the 
section~\ref{bntifent}. The notations are fixed in the section~\ref{cadre}
and our main result is stated in the section~\ref{princ}.

%--------------------------------------------------------------------------
\section{Setting}\label{cadre}
%--------------------------------------------------------------------------
%--------------------------------------------------------------------------

Let us briefly fix the notations (for details see~\cite{mout3}).
We consider a \textit{tree} $(S,A)$
i.e. a non-oriented, locally finite, connected and simply
connected graph with \textit{vertices} in $S$ and \textit{edges} in $A$.
We will use the usual notions of \textit{path}, \textit{distance} and
\textit{geodesic path} and note $x\sim y$ iff $(x,y)\in A$.

We also consider a transient random walk $(X_n)_n$ on $S$ such that the 
transition probability $p(x,y)>0$ iff $x\sim y$. Denote 
by $P_x$ the distribution of the walk starting from $x$ and by $p_n(x,y)$
the probability $P_x[X_n=y]$ of reaching $y$ from $x$ in $n$ steps. 

The \textit{Green function} 
$G(x,y)= \sum_{n=0}^\infty p_n(x,y)$ is finite by transience. 
Denote by $H(x,y)$ the probability of reaching $y$ starting from $x$. 
If $z$ is on the geodesic path $[x,y]$, the simple connectivity implies
\begin{equation}\label{formg}
  H(x,y)=H(x,z)H(z,y)  
  \ \ \ \mathrm{and}\ \ \ 
  G(x,y)=H(x,z)G(z,y).  
\end{equation}
If $U\subset S$, the \textit{Green function of} $U$, defined on
$U\times U$, is the expectation of the number
of times the walk starting from $x$ hits $y$ before exiting $U$. 

The \textit{Laplacian} of a function $f$ on $S$ is
$\Delta f(x)= E_x[f(X_1)]-f(x)$. 
The function $f$ is \textit{harmonic} if $\Delta f = 0$.

Let $u$ be a fixed harmonic function.
The \textit{stochastic energy} of $u$ is 
$    J^*(u) = \sum_{k=0}^\infty \left( \Delta u^2\right)(X_k) $
(non-negative terms).
The events $\calL^{**}$, $\calN^{**}$ and $\calJ^{**}$
are defined respectively by the convergence of $(u(X_n))_n$,
its boundedness and the finiteness of the stochastic energy.
The Martingale theorem implies $ \calJ^{**} \simsubset \calL^{**}$
($P_x$-almost sure inclusion)~\cite{mout3}. 
It is known since P.~Cartier~\cite{cart} that geometric and 
Martin compactifications agree and
the random walk converges almost surely to a point of the boundary
$\dS$.
The exit law starting from $x$ is the harmonic measure $\mu_x$ and 
$\mu=(\mu_x)_x$ is a familly of equivalent measures.
Conditioning by Doob's method of $h$-processes gives 
probabilities $P_x^{\theta}$ (ending at $\theta$).
Asymptotic events verify $0$--$1$ law and we define sets
$  \calL^{*} 
      = \left\{
          \theta \in \dS |P_x^\theta(\calL^{**})=1
        \right\}                   
$, 
$    \calN^{*} 
      = \left\{
          \theta \in \dS | P_x^\theta(\calN^{**})=1
        \right\}                   
$,
$    \calJ^{*} 
      = \left\{
          \theta \in \dS | P_x^\theta(\calJ^{**})=1
        \right\}                   
$,
which determine stochastic notions of convergence, boundedness and 
finiteness of the energy at $\theta\in\dS$. For $\theta \in \calL^{*}$, 
$\lim u(X_n)$ is $P_x^\theta$-a.s. constant (independent from $x$) 
and called the \textit{stochastic limit} at $\theta$.

Fix a base point $o$.
For $\theta\in\dS$, $\gamma_\theta$ is the geodesic ray from $o$ 
to $\theta$ and for $c\in\N$,
$\Gamma_c^{\theta}=\{y\in S | d(y,\gamma_{\theta})\leq c\}$
is a \textit{non-tangential tube}.
Let $u$ be a harmonic function.
For $c\in\N$, its \textit{$c$-non-tangential energy} 
at $\theta$ is 
$ J_c^\theta(u)= \sum_{y\in\Gamma_c^{\theta}} \Delta u^2 (y)$ 
and its \textit{radial energy} at $\theta$ is 
$J^\theta(u)= J_0^\theta(u)=\sum_{k=0}^\infty \Delta u^2 (\gamma_\theta(k)) $.
There is \textit{radial} convergence,
boundedness or finiteness of the energy depending wether
$(u(\gamma_\theta(n)))_n$ converges, is bounded or
has finite radial energy.  
There is \textit{non-tangential} convergence
of $u$ at $\theta$ if for all $c\in\N$, $u(y)$ has a limit
when $y$ goes to $\theta$ staying in $\Gamma_c^{\theta}$.
There is \textit{non-tangential} boundedness (resp. 
finiteness of the energy) if 
for all $c\in\N$, $u$ is bounded on $\Gamma_c^{\theta}$
(resp. $J_c^\theta(u) < +\infty$).

%--------------------------------------------------------------------------
\section{Main result}\label{princ}
%--------------------------------------------------------------------------
%--------------------------------------------------------------------------

We now suppose $\Hyp$:
$    \exists \varepsilon > 0, 
     \exists \eta > 0,
     \forall x\sim y, 
     \varepsilon \leq p(x,y) \leq \frac{1}{2} - \eta     $,
a discrete analogue of the pinched curvature for manifolds.
It also forces at least three neighbors for each vertex, 
and ensures transience.
We proved in~\cite{mout3}:
\begin{theorem}
  For a harmonic function $u$ on a tree with random walk satisfying $\Hyp$,
  the notions of radial convergence, radial boundedness, radial 
  finiteness of the energy, 
  stochastic convergence, stochastic boundedness, stochastic finiteness of the 
  energy, are $\mu$-almost equivalent. 
\end{theorem}
We prove here the following theorem:
\begin{theorem}
  Under the same hypotheses, the notions of non-tangential convergence,
  non-tangential boundedness and non-tangential finiteness of the energy 
  are $\mu$-almost equivalent to the notions above.
\end{theorem}

Considering the trivial implications, it is sufficient to prove 
that stochastic convergence implies non-tangential convergence and 
non-tangential boundedness implies almost surely non-tangential finiteness
of the energy.

%--------------------------------------------------------------------------
\section{Stochastic implies NT convergence}\label{csicnt} 
%--------------------------------------------------------------------------
%--------------------------------------------------------------------------

The first implication needs 
the following lemma due to A.~Ancona in a general setting \cite{anco3}, 
but easily proved here by simple connectivity:
\begin{lemma}
  If $(x_n)_n$ is a sequence converging non-tangentially to 
  $\theta\in\dS$, the walk hits $P_o^\theta$-a.s. infinitely many $x_n$.
\end{lemma}

Let us see how this lemma helps. Assume that the harmonic function $u$
has a stochastic limit $l\in\R$ at $\theta$ but does not converge 
non-tangentially towards $l$ at $\theta$. There exists  $\delta>0$
and a sequence $(x_n)_n$ converging non-tangentially to 
$\theta$ such that $|u(x_n)-l|\geq \delta$ for all $n$.
As the random walk  $(X_k)_k$ hits $P_o^\theta$-a.s. infinitely many
$x_n$ by the lemma, one can extract a subsequence $(X_{k_j})_j$ such that
$|u(X_{k_j})-l|\geq \delta$ for all $j$. 
Hence, $P_o^\theta$-almost surely, the function $u$ does not converge
towards $l$ along $(X_k)_k$ which leads to a contradiction.

Le us now prove the lemma. Recall that the principle of the method of Doob's 
$h$-processes is to consider a new Markov chain defined by
$p^\theta(x,y)=\frac{K_\theta(y)}{K_\theta(x)}p(x,y)$ where
the Martin kernel $K_\theta(x)$ is defined as
$\lim_{y \rightarrow\theta}\frac{G(x,y)}{G(o,y)}$
(see for example~\cite{dynk}).
This formula leads to analogous fomulae for the 
 $p_n^\theta$ and the associated functions $H^\theta$ and $G^\theta$.
Consider for a fixed $n$ the projection $y_n$ of $x_n$ on
the geodesic ray  $\gamma_\theta$ (see~\cite{mout3}).
As the random walk starting from $o$ and conditioned to end at $\theta$
hits almost surely $y_n$ due to the tree structure, the strong Markov property
gives 
$ H^\theta(o,x_n)= H^\theta(y_n,x_n)= 
   \frac{K_\theta(x_n)}{K_\theta(y_n)}H(y_n,x_n) $.
By definition of the Martin kernel, 
$ \frac{K_\theta(x_n)}{K_\theta(y_n)} 
   = \lim_{y \rightarrow\theta}\frac{G(x_n,y)}{G(y_n,y)} $
and $G(x_n,y)=H(x_n,y_n)G(y_n,y)$ as soon as $y_n\in[x_n,y]$, so
$H^\theta(o,x_n)= H(x_n,y_n)H(y_n,x_n)$.
The distance between $x_n$ and $y_n$ is bounded as $(x_n)_n$
converges non-tangentially to $\theta$, hence the last 
product is bounded from below by a constant  $C>0$ using $\Hyp$.
By Fatou's lemma, the probability conditioned to end at $\theta$
of hitting infinitely many $x_n$ is not smaller than $C$ and
the asymptotic $0$-$1$ law ensures that it equals $1$, which completes 
the lemma's proof.

%--------------------------------------------------------------------------
\section{NT boundedness implies finite NT energy}\label{bntifent}
%--------------------------------------------------------------------------
%--------------------------------------------------------------------------

Denoting $\calN_c=\{\theta \in\dS | \sup_{\Gamma_c^\theta} |u| < +\infty\}$
and $\calJ_c=\{\theta \in\dS | J_c^\theta(u)< +\infty\}$, 
we will show that for all $c\in\N$, $\calN_{c+1}\simsubset\calJ_c$,
which will give the wanted result by monotonous intersection. 
Let us write
$\calN_{c+1}= \bigcup_{N\in \N} \calN_{c+1}^N$, where
$$ \calN_{c+1}^N =        \left\{
                           \theta\in\dS
                           \left|\ 
                             \sup_{\Gamma_{c+1}^\theta} |u| \leq N
                           \right.  
                         \right\}.$$   
By countability it is sufficient to prove that for all 
$N$, $\calN_{c+1}^N\simsubset\calJ_c$.
Let us fix $N\in\N$. Denote 
$\Gamma=\bigcup_{\theta\in\calN_{c+1}^N}\Gamma_{c}^\theta $ and $\tau$
the exit time from $\Gamma$. As 
$$ M_n = u^2(X_n)-\sum_{k=0}^{n-1} \Delta u^2(X_k)  $$
is a martingale (see~\cite{mout3}), 
Doob's stopping time theorem for the bounded exit time
$\tau \wedge n$ gives
$ E_o\left[M_{\tau\wedge n}\right] 
      = E_o[M_0] 
      = u^2(o) \geq 0  $,
hence
$$ E_o\left[\sum_{k=0}^{{\tau\wedge n}-1} \Delta u^2(X_k)\right] 
\leq E_o\left[u^2(X_{\tau\wedge n})\right]
.$$
As $X_{\tau\wedge n}$ is at distance at most $1$ from $\Gamma$,
it lies in a tube $\Gamma_{c+1}^\theta$ where 
$\theta\in\calN_{c+1}^N$ and $|u(X_{\tau\wedge n})| \leq N$.
When $n$ goes to $\infty$, monotonous convergence
($\Delta u^2 \geq 0$) and the desintegration formula (see~\cite{mout3})
give then, for $\mu$-almost all $\theta\in\dS$,
$$ E_o^\theta\left[\sum_{k=0}^{\tau-1} \Delta u^2(X_k)\right] <+\infty  
.$$
Let us use a conditioned version of formula~2 from~\cite{mout3},
which will be proved later~:
\begin{lemma}\label{forproanacond}
For a function $\varphi\geq 0$ on $\Gamma$ and
$\tau$ the exit time of $\Gamma$,
$$   E_o^\theta\left[\sum_{k=0}^{\tau-1} \varphi(X_k)\right]
     = \sum_{y\in \Gamma} \varphi(y) G_\Gamma(o,y)K_\theta(y).  $$  
\end{lemma}
This lemma implies that for $\mu$-almost all $\theta\in\dS$,
$\sum_{y\in \Gamma} \Delta u^2(y) G_\Gamma(o,y)K_\theta(y)$ is finite.
In order to get an energy, we will show that
$G_\Gamma(o,y)K_\theta(y)$ 
is bounded from below using the two following lemmas.
The first one is due to A.~Ancona~\cite{anco3} but has
a very simple proof in the present context of trees.
The second one enables comparison between $G_\Gamma$ and $G$.
\begin{lemma}\label{lemanco}
$\forall c\in\N,\exists\alpha>0,\forall\theta\in\dS,
\forall y\in\Gamma_c^\theta,G(o,y)K_\theta(y)\geq \alpha $.
\end{lemma}
\begin{lemma}\label{lemcvd}
  For $U \subset S$ containing $\Gamma_c^\theta$ and $\tau$ 
  the exit time of $U$,
  $$ \lim_{y\in\Gamma_c^\theta, y\rightarrow\theta} 
         \frac{G_U(o,y)}{G(o,y)} = P_o^\theta[\tau=+\infty] .$$
\end{lemma} 
By lemma~\ref{lemanco}, for $\mu$-almost all $\theta\in\calN_{c+1}^N$,
$$  \sum_{y\in \Gamma_c^\theta} 
      \Delta u^2(y) \frac{G_\Gamma(o,y)}{G(o,y)}  <+\infty  .$$
If we show that for $\mu$-almost all $\theta\in\calN_{c+1}^N$,
$P_o^\theta[\tau=+\infty]>0$, lemma~\ref{lemcvd} gives  
$\calN_{c+1}^N\simsubset\calJ_c$. 
The proof of that fact is the same as in the analogous 
radial proof~\cite{mout3} which completes the theorem's proof.

Let us now prove the lemmas. Concerning lemma~\ref{forproanacond}, using Fubini,
$$   E_o^\theta\left[\sum_{k=0}^{\tau-1} \varphi(X_k)\right]
     = \sum_{k=0}^{\infty} E_o^\theta\left[ \varphi(X_k)\un{k<\tau}\right]
.$$
The random variable $\varphi(X_k)\un{k<\tau}$ being measurable with respect to 
the $\sigma$-algebra generated by $(X_i)_{i\leq k}$ 
(see~\cite{mout3}) and using formula~2 from~\cite{mout3}, the
expectation above equals
$$  \sum_{k=0}^{\infty} E_o\left[ \varphi(X_k)\un{k<\tau} K_\theta(X_k)\right] 
  = E_o\left[\sum_{k=0}^{\infty} \varphi(X_k)\un{k<\tau} K_\theta(X_k)\right] $$
$$  = \sum_{y\in \Gamma} \varphi(y) G_\Gamma(o,y)K_\theta(y) ,$$
which finishes the proof of lemma~\ref{forproanacond}.

Let us prove lemma~\ref{lemanco}.
Denote $\pi(y)$
the projection of $y$ on $\gamma_\theta$ (see~\cite{mout3})
and remark that for
$z\in (\pi(y),\theta)$, $G(o,z)=H(o,\pi(y))G(\pi(y),z)$ and
$G(y,z)=H(y,\pi(y))G(\pi(y),z)$ by formula~\ref{formg}.
Hence $\frac{G(y,z)}{G(o,z)}=\frac{H(y,\pi(y))}{H(o,\pi(y))}$ 
does not depend anymore on $z$ and its limit when
$z$ goes to $\theta$ is then
$K_\theta(y)=\frac{H(y,\pi(y))}{H(o,\pi(y))}$. By formula~\ref{formg},
$$G(o,y)K_\theta(y)=H(y,\pi(y))\frac{G(o,y)}{H(o,\pi(y))}
=H(y,\pi(y))H(\pi(y),y)G(y,y).$$ 
But $G(y,y)\geq p_2(y,y)\geq 3\varepsilon^2$ and
$H(y,\pi(y))H(\pi(y),y)\geq \varepsilon^{2c}$
by $\Hyp$ and $d(y,\pi(y))\leq c$, which finishes the 
proof of lemma~\ref{lemanco}.

Let us prove lemma~\ref{lemcvd}~:
$$ G_U(o,y)= G(o,y)-E_o[G(X_\tau,y)\un{\tau<+\infty}] $$
$$= G(o,y)\left(1-
              E_o\left[\frac{G(X_\tau,y)}{G(o,y)}\un{\tau<+\infty}
                 \right]
        \right)  $$
and by definition of Martin's kernel, if we could switch
the limit and expectation, by a conditioning formula~\cite{mout3},
$$ \lim_{y\in\Gamma_c^\theta, y\rightarrow\theta} 
         \frac{G_U(o,y)}{G(o,y)} 
= 1- E_o[K_\theta(X_\tau)\un{\tau<+\infty}]  
 = P_o^\theta[\tau=+\infty] .$$
We now justify that inversion by Lebesgue's theorem.
The idea is to bound, when $\tau$ is finite, 
$\frac{G(X_\tau,y)}{G(o,y)}$ by a multiple of
$K_\theta(X_\tau)$. We compare for that purpose
$G(X_ \tau,y)$ with $K_\theta(X_\tau)$.   
Denote again by $\pi$ the projection function on $\gamma_\theta$.
We distinguish two cases 

If $\pi(X_\tau)\in [o,\pi(y)]$,
$ \frac{G(X_\tau,y)}{K_\theta(X\tau)}
  = \frac{G(\pi(X_\tau),y)}{K_\theta(\pi(X\tau))}
  = \frac{G(o,y)}{K_\theta(o)} = G(o,y) $,
by formula~\ref{formg} and the remark that this formula also implies
by definition of $K_\theta$ and by taking the limit that
$K_\theta(X_\tau)=H(X_\tau,\pi(X_\tau))K_\theta(\pi(X_\tau))$ 
and $K_\theta(o)=H(o,\pi(X_\tau))K_\theta(\pi(X_\tau))$. 

If $\pi(X_\tau)\not\in [o,\pi(y)]$,
again
$ \frac{G(X_\tau,y)}{K_\theta(X\tau)}
= \frac{G(\pi(X_\tau),y)}{K_\theta(\pi(X\tau))}$.
We also have, by definition and formula~\ref{formg},
$K_\theta(\pi(X_\tau))=(H(o,\pi(X_\tau)))^{-1}$, 
hence the quotient above equals
$ H(o,\pi(X_\tau))G(\pi(X_\tau),y)
  = H(o,\pi(y))H(\pi(y),\pi(X_\tau))G(\pi(X_\tau),y)$.
We know that $G$ is bounded (see~\cite{pica-woes,mout3})
and $H$ is a probability, so it just remains to compare
$H(o,\pi(y))$ with $G(o,y)$. But
$\frac{H(o,\pi(y))}{G(o,y)}= (G(\pi(y),y))^{-1}$ and $\frac{1}{G}$ is 
bounded by $\frac{1}{3\varepsilon^2}$.

Merging the two cases gives a constant $\beta$
such that $\frac{G(X_\tau,y)}{K_\theta(X\tau)} \leq \beta G(o,y)$,
which enables to use Lebesgue's theorem and completes the proof of 
lemma~\ref{lemcvd}.

% Les remerciements sont dans une section, sans numérotation
%\section*{Remerciements}
% Remerciements - texte ici

\vspace{5mm}
\noindent
{\small 
Fr\'ed\'eric Mouton\\ 
Universit\'e Grenoble~1\\ 
BP 74\\
38402 Saint-Martin-d'H\`eres Cedex (France)\\
\texttt{Frederic.Mouton@ujf-grenoble.fr}
}


\begin{thebibliography}{00}

\bibitem{anco3}
A.~Ancona.
\newblock Th{\'e}orie du potentiel sur les graphes et les vari{\'e}t{\'e}s.
\newblock In P.L. Hennequin, editor, {\em {\'E}cole d'{\'e}t{\'e} de
  probabilit{\'e}s de {S}aint-{F}lour {XVIII}}. Springer Lect. Notes in Math.
  1427, Berlin, 1990.

\bibitem{atan-pica}
L.~Atanasi and M.A.~Picardello.
\newblock The lusin area function and local admissible convergence of harmonic
  functions on homogeneous trees.
\newblock {\em Trans. of A.M.S.}, 360:3327--3343, 2008.

\bibitem{cart}
P.~Cartier.
\newblock Fonctions harmoniques sur un arbre.
\newblock In {\em Symposia Mathematica}, volume~IX, pages 203--270. Academic
  Press, London and New-York, 1972.

\bibitem{dynk}
E.B. Dynkin.
\newblock Boundary theory of markov processes (the discrete case).
\newblock {\em Russ. Math. Surv.}, 24:1--42, 1969.

\bibitem{mout2}
F.~Mouton.
\newblock Comportement asymptotique des fonctions harmoniques en courbure
  n\'egative.
\newblock {\em Comment. Math. Helvetici}, 70:475--505, 1995.

\bibitem{mout3}
F.~Mouton.
\newblock Comportement asymptotique des fonctions harmoniques sur les arbres.
\newblock {\em S\'eminaire de Probabilit\'es, Universit\'e de Strasbourg},
  XXXIV:353--373, 2000.

\bibitem{pica-woes}
M.~A. Picardello and W.~Woess.
\newblock Finite truncations of random walks on trees.
\newblock In {\em Symposia Mathematica}, volume XXIX, pages 255--265. Academic
  Press, London and New-York, 1987.

\end{thebibliography}
\end{document}